\documentclass[12pt]{article}
\usepackage{amscd}
\newtheorem{Definition}{Definition}
\newtheorem{Proposition}{Proposition}
\newtheorem{Lemma}{Lemma}

\newtheorem{Example}{Example}
\newtheorem{Theorem}{Theorem}

\newenvironment{keywords}{
       \list{}{\advance\topsep by0.35cm\relax\small
       \leftmargin=1cm
       \labelwidth=0.35cm
       \listparindent=0.35cm
       \itemindent\listparindent
       \rightmargin\leftmargin}\item[\hskip\labelsep
                                     \bfseries Keywords:]}
     {\endlist}

\begin{document}

\title{IMPRIMITIVITY THEOREM FOR GROUPOID REPRESENTATIONS}

\author{
Leszek Pysiak \\
Department of Mathematics and Information Science, \\ Warsaw
University of Technology \\
Plac Politechniki 1, 00-661 Warsaw, Poland}

\date{March 30, 2010}
\maketitle

\begin{abstract}
We define and investigate the concept of the groupoid representation induced by a representation of the isotropy subgroupoid. Groupoids in question are locally compact transitive topological groupoids. We formulate and prove the imprimitivity theorem for such representations which is a generalization of the classical Mackey's theorem known from the theory of group representations.

% \PACS{98.80.-k \and 98.80.Bp \and 98.80.Qc}
% \subclass{MSC code1 \and MSC code2 \and more}

\end{abstract}

\section{Introduction}
The present paper, devoted to the study of the theory of groupoid representations, is a continuation of my previous work \cite{Pysiak} in which one can find a presentation of the groupoid concept and of the groupoid representation concept, important examples as well as relationships between groupoid representations and induced group representations (see also \cite{Paterson}, \cite{Renault}, \cite{Weinstein}, \cite{Landsman}). Groupoids have now found a permanent place in manifold domains of mathematics, such as: algebra, differential geometry, in particular noncommutative geometry, and algebraic topology, and also in numerous applications, notably in physics. It is a natural tool to deal with symmetries of more complex natura than those described by groups (see \cite{Weinstein}, \cite{Landsman}, \cite{Brown}). Groupoid representations were investigated by many authors and in many ways (see \cite{Westman}, \cite{Renault}, \cite{Paterson}, \cite{Bos}, \cite{Buneci}).
\par
In a series of works (\cite{Full}, \cite{Conceptual},
\cite{Random}, \cite{PysiakTime}) we have developed a model
unifying gravity theory with quantum mechanics in which it is a
groupoid that describes symmetries of the model, namely the
transformation groupoid of the pricipal bundle of Lorentz frames
over the spacetime. To construct the quantum sector of the model
we have used a regular representation of a noncommutative
convolutive algebra on this groupoid in the bundle of Hilbert
spaces. In paper (\cite{Malicious}) we have applaied this groupoid
representation to investigate spacetime singulariies, and in
\cite{Aharonov} the representation of the fundamental groupoid to
the gravitational Aharonov-Bohm effect.

\par
The present work is aimed at introducing the concept of the
groupoid representation induced by a representation of the
isotropy subgroupoid. We assume that the groupoid in question is a
locally trivial topological groupoid and as a topological space it
is a locally compact Husdorff space (Section 4). This concept,
framed ``in the spirit of Mackey'' is a natural generalization of
induced represetation of locally compact groups, created and
investigated by him \cite {Mackey3}. Representations, investigated
in the present work, are unitary and are realized in Hilbert
bundles over the unit spaces of a given groupoid \cite{Paterson}.
Section 5 contains the formulation and the proof of the
imprimitivity theorem  for groupoids which is a generalization of
the classical Mackey's imprimitvity theorem for group
representations \cite{Mackey2}, \cite{Mackey3}. The theorem says
that every unitary groupoid representation, for which there exists
the imprimitivity system, is a representation induced by a
representation of the isotropy subgroupoid.
\par
In Section 6, I investigate induced representations of the
transformation gropoid over a homogeneous space of the group $G$
and show that there exists a strict connection between these
representations and induced representations of the group $G$ (in
the sense of Mackey.
\par
In Section 7, I give a physical interpretation of concepts
analyzed in Section 6. I describthe a representation that has been
used in the mentioned above model unifying gravity and quanta when
this model is reduced (as the result of the act of measurement) to
the usual quantum mechanics. And then I consider the
energy-momentum space for a massive particle (it is a homogeneous
space of the Lorentz group) and the transformation groupoid
corresponding to this space.  I also give a definition
 (in the sense of Mackey \cite{Mackey2}, \cite{Mackey1}) of a particle as
an imprimitivity sestem for the unitary representation of this groupoid.

\section{Preliminaries}
Let ${\cal G}$ be a groupoid over a set $X$ (the base of ${\cal
G})$. We recall (cf. \cite{SilvaWeinstein}, \cite{Paterson} ) that a
groupoid ${\cal G}$ is a set with a partially defined multiplication
"$\circ$" on a subset ${\cal G}^2$ of ${\cal G} \times {\cal G}$,
and an inverse map $g \rightarrow g^{-1}$ defined for every $g \in
{\cal G}$. The multiplication is associative when defined. One has
an embedding $\epsilon: X \rightarrow {\cal G}$ called the identity
section and two structure maps $ d,r: {\cal G}\rightarrow X$ such
that
  \[\epsilon(d(g))=g^{-1}\circ g\]
\[\epsilon(r(g))=g\circ g^{-1}\]
for $g \in{\cal G}.$
\par
Let us introduce the following fibrations in the set ${\cal G}$:
\[{\cal G}_x= \{g\in {\cal G}: d(g)=x\} \]
\[{\cal G}^x = \{g\in {\cal G}: r(g)=x\} \] for $x \in X$. Let us
also denote ${\cal G}_x^y = {\cal G}^x \bigcap {\cal G}_y $, and
consider the set ${\cal G}_x^x = {\cal G}^x \bigcap {\cal
G}_x $ for $x\in X$. It has the group structure and is called
\emph{the isotropy group of the point $x$}. It is clear that the
set $\Gamma = \bigcup_{x \in X} {\cal G}_x^x$ has the structure of
a subgroupoid of ${\cal G}$ over the base $X$ (all the structure
maps are the restrictions of the structure maps of
${\cal G}$ to $\Gamma$).
\par
We call ${\cal G}$ \emph{a transitive groupoid}, if for each pair
of elements $x_1, x_2 \in X$ there exists $g \in{\cal G}$ such
that $d(g)= x_1$ and $r(g) = x_2$.
\par
A groupoid ${\cal G}$ is \emph{a topological groupoid} if ${\cal
G}$ and $X$ are topological spaces and all structure maps are
continuous (in particular, the embedding $\epsilon$ is a
homeomorphism of $X$ onto its image).
\par
In the following we assume that ${\cal G}$ (and thus $X$) is a
locally compact Hausdorff space.
\begin{Example}
\emph{ A pair groupoid}. Let $X$ be a locally compact Hausdorff
space. Take ${\cal G}= X\times X$. We define the set ${\cal G}^2$
of composable elements as  ${\cal G}^2= \{((x,y),(y,z)): x,y,z \in
X\}\subset{\cal G}\times {\cal G}$  and a multiplication, for
$((x,y),(y,z))\in {\cal G}^2 $, by
\[ (x,y)\circ(y,z)= (x,z).\]
Moreover, we have: $(x,y)^{-1}=(y,x)$, $d(x,y)=y$, $r(x,y)=x$,
$\epsilon(x)=(x,x)$. With such defined structure maps ${\cal G}$
is a groupoid, called \emph{pair groupoid}.
\end{Example}
\begin{Example}
\emph{A transformation groupoid.} Let $X$ be a locally compact
Hausdorff space, and $G$ a locally compact group. Let $G$ act
continuously on $X$ to the right, $X\times G\rightarrow X$. Denote
$(x,g)\mapsto xg $. We introduce the groupoid structure on the set
${\cal G}= X\times G$ by defining the following structure maps.
The set of composable elements ${\cal G}^2= \{((xg,h),(x,g): x\in
X, g,h\in G \}\subset{\cal G}\times {\cal G}$, and the
multiplication for $((xg,h),(x,g))\in {\cal G}^2 $ is given by
\[(xg,h)\circ(x,g)= (x,gh).\]
And also $(x,g)^{-1}=(xg,g^{-1})$, $d(x,g)=x$, $r(x,g)=xg$,
$\epsilon(x)=(x,e_G)$. This groupoid is called \emph{the
transformation groupoid.}
\end{Example}
\par
Let us recall (cf.\cite{Paterson}) the concept of \emph{right Haar
System.}
\par
\begin{Definition}
\emph{A right Haar system} for the groupoid ${\cal G}$ is a family
$\{\lambda_x\}_{x \in X}$ of regular Borel measures defined on the
sets ${\cal G}_x$ (which are locally compact Hausdorff spaces)
such that the following three conditions are satisfied:
\begin{enumerate}
\item
the support of each $\lambda_x$ is the set ${\cal G}_x$,
\item
(continuity) for any $f \in C_c({\cal G})$ the function $f^0$,
where
\[ f^0(x)=\int_{{\cal G}_x} f d\lambda_x ,\]
belongs to $C_c(X)$,
\item
(right invariance) for any $g \in {\cal G}$ and any $f \in
C_c({\cal G})$,
\[ \int_{{\cal G}_{r(g)}} f(h \circ g)
d\lambda_{r(g)}(h)=\int_{{\cal G}_{d(g)}} f(u) d\lambda_{d(g)}(u).
\]

\end{enumerate}
\end{Definition}
One can also consider the family $\{\lambda^x\}_{x \in X}$ of
left-invariant measures, each $\lambda^x$ being defined on the set
${\cal G}^x$ by the formula $\lambda^x(E)=\lambda_x(E^{-1})$ for
any Borel subset $E$ of ${\cal G}^x$ (where $E^{-1}=\{g \in {\cal
G}: g^{-1} \in E\}$). Then the invariance condition assumes the
form:
\[ \int_{{\cal G}^{d(g)}}f(g\circ h)d\lambda^{d(g)}(h) = \int_{{\cal
G}^{r(g)}}f(u)d\lambda^{r(g)}(u). \]
 Now, let $\mu$ be a regular
Borel measure on $X$. We can consider the following measures which
will be called \emph{measures associated with $\mu$}: $\nu=\int
\lambda_x d\mu(x)$ on ${\cal G}$, $\nu^{-1}=\int \lambda^x
d\mu(x)$ and $\nu^2=\int \lambda_x \times \lambda^x d\mu(x)$ on
${\cal G}^2$.
\par
If $\nu = \nu^{-1}$ we say that the measure $\mu$ is \emph{a
${\cal G}$-invariant measure on $X$.}
\par
\begin{Definition}
A topological groupoid ${\cal G}$ on $X$ is called \emph{locally
trivial} if there exist a point $x \in X$, an open cover $\{U_i\}$
of $X$ and continuous maps $s_{x,i}: U_i \to {\cal G}_x$ such that
$r \circ s_i = id_{U_i}$ for all $i$.
\end{Definition}

\begin{Proposition} Assume that ${\cal G}$ is a locally trivial
groupoid on $X$ and $X$ is second countable space. Let $\mu$ be a
regular Borel measure on $X$. Then
\begin{enumerate}
\item
${\cal G}$ is transitive,
\item
all isotropy groups of ${\cal G}$ are isomorphic with each other,
\item
for every $y \in X$ there exist an open cover $\{V_j\}$ of $X$
and continuous maps $s_{y,j}: V_j \to {\cal G}_y$ such that $r
\circ s_j = id_{v_j}$,
\item
for every $x \in X$ there exists a section $s_x: X \to {\cal G}_x$
which is $\mu$-measurable, i.e., for every Borel set $B$ in ${\cal
G}_x$, $s_x^{-1}(B)$ is $\mu$-measurable subset of $X$,
\item
if the measure $\mu$ has the property that $\mu (\overline{A})=
\mu(A)$ for every $\mu$-measurable subset $A$ of $X$, then the
section $s_x$ is $\mu-a.e.$ continuous on $X$.
\end{enumerate}
\end{Proposition}
\textbf{Proof:}
\begin{enumerate}
\item
Let $y_1, y_2 \in X$. Suppose that $y_1 \in U_1$ and $y_2 \in
U_2$. Then $r(s_1(y_1))=y_1$ and $r(s_2(y_2))=y_2$. But
$g=s_2(y_2) \circ s_1(y_1)^{-1}$ has the property $d(g)=y_1$ and
$r(g)=y_2$, and this means that ${\cal G}$ is transitive.
\item
For $x,y\in X $ let $g_{yx}$ be an element of ${\cal G}$ such that
$d(g_{yx})=x$ and $r(g_{yx})=y$. Then we have an isomorphism of
the isotropy groups ${\cal G}_x^x$ and ${\cal G}_y^y$ given by the
formula ${\cal G}_x^x \ni \gamma \mapsto g_{yx}\circ\gamma\circ
g_{yx}^{-1}\in {\cal G}_y^y$.
\item
Let $g_{xy}$ be an element of ${\cal G}$ such that $d(g_{xy})=y$
and $r(g_{xy})=x$. Then in the fiber ${\cal G}_y$ we can simply
define $s_{y,i}(z)=s_{x,i}(z) \circ g_{x,y}$ for $z \in U_i$.
\item
Since $X$ is second countable space, we can take a countable cover
$\{U_i\}_{i=1,2,...}$ of the space $X$. Let us define $s_x(z) =
s_1(z)$ for $z \in U_1$, $s_x(z) = s_2(z)$ for $z \in U_2
\setminus U_1 , \ldots , s_x (z) = s_n (z)$ for $z \in U_n
\setminus (U_1 \cup U_2 \cup \ldots U_{n-1})$ etc.
\par
It is easily seen that $s_x(z)$ is measurable.
\item
The set of discontinuity of $s_x$ is contained in the union of
sets $\partial U_i = \overline{U_i}\setminus U_i$, $i=1,2,...$,
which is of measure zero.
\par
This ends the proof.
\end{enumerate}
\par
From now on we assume that considered groupoids are locally
compact and satisfy the assumptions of Proposition 1. It is known
that in the case of any locally trivial groupoid there exists a
right Haar system (see \cite{Paterson}). Let us choose a
collection of sections $\{s_x\}_{x\in X}$, defined by Proposition
1.3, and denote by $d\gamma_y$ the right Haar measure on the
isotropy group $\Gamma_y$. (From the assumption that the groupoid
${\cal G}$ is locally compact it follows that all isotropy groups
are locally compact and have Haar measures.)

\begin{Definition}
A right Haar system $\{\lambda_x\}_{x\in X}$ on the groupoid
${\cal G}$ is called \emph{consistent with  a Borel regular
measure $\mu$ on the base space $X$}  if, for every $x\in X$ and
any $f \in C_c({\cal G}_x)$,
\[ \int_{{\cal G}_x} f(g)
d\lambda_x(g)=\int_X\int_{\Gamma_y} f(\gamma \circ s_x(y))
d\gamma_y d\mu(y).
\]
\end{Definition}
The above formula gives an explicit construction of the right Haar
system for many classes of groupoids (see, Section 3 for pair
groupoid, and Section 6 for transformation groupoid).
\par
Let us recall the concept of \emph{groupoid representation}
\cite{Paterson},\cite{Renault}. It involves a Hilbert bundle
$\mathbf{H}$ over $X$, $\mathbf{H}=(X, \{H_x\}_{x \in X}, \mu)$
(Dixmier \cite{Dixmier} uses for it the name of $\mu$-measurable
field of Hilbert spaces over $X$). Here all Hilbert spaces ${H_x}$
are assumed to be separable.
\par
Let $\mu$ be a ${\cal G}$-invariant measure on $X$, and $\nu$ and
$\nu^2$ the associated measures on ${\cal G}$.
\begin{Definition}
\emph{A unitary representation of a groupoid ${\cal G}$} is the
pair $({\cal U},\mathbf{H})$ where $\mathbf{H}$ is a Hilbert
bundle over $X$ and \ ${\cal U}=\{U(g)\}_{g \in {\cal G}}$ is a
family of unitary maps \ $U(g): H_{d(g)} \rightarrow H_{r(g)}$
such that:
\begin{enumerate}
\item
$U(\epsilon(x)) = id_{H_x}$ for all $x \in X$,
\item
$U(g)\circ U(h)=U(g\circ h)$ for $\nu^2 - a.e. \ (g,h) \in {\cal
G}^2$,
\item
$U(g^{-1})=U(g)^{-1} \ for \ \nu - a.e. \ g \in {\cal G}$,
\item
For every $\phi, \psi \in L^2(X, \mathbf{H}, \mu)$,
\[ {\cal G} \ni g \rightarrow (U(g) \phi(d(g)), \psi(r(g)))_{r(g)} \in
\mathcal{C}\] is $\nu$-measurable on ${\cal G}$. (Here $L^2(X,
\mathbf{H}, \mu)$ denotes the space of square-integrable sections
of the bundle $\mathbf{H}$, and \ $(\cdot,\cdot)_x$ denotes the
scalar product in the Hilbert space $H_x$.)
\end{enumerate}
\end{Definition}

\section{Elementary properties of representations of groupoids.}

\begin{Definition}
Unitary representations $ ({\cal U}_1,\mathbf{H}_1)$ and $({\cal
U}_2,\mathbf{H}_2)$ of a groupoid $\cal G$ are said to be
\emph{unitarily equivalent} if there exists a family
$\{A_x\}_{x\in X}$ of isomorphisms of Hilbert spaces $A_x:
H_{1x}\rightarrow H_{2x}, x\in X$ such that for every $x,y\in X$
and for $\nu-$a.e. $g\in {\cal G}_x^y$ the following diagram
commutes
\[
\begin{CD}
H_{1x} @>{U_1(g)}>> H_{1y} \\
@V{A_x}VV @VV{A_y}V \\
H_{2x} @>>{U_2(g)}> H_{2y}
\end{CD}
\]
\end{Definition}
\begin{Definition}
Let $ ({\cal U},\mathbf{H})$ be an unitary representation of $\cal
G$ and let $\mathbf{H}_1$ be a Hilbert subbundle of \
$\mathbf{H}$. We say that $\mathbf{H}_1$ is $\cal G$-invariant if
$U(g) H_{1,x} \subset H_{1,y}$ for every $x,y\in X$ and for
$\nu-$a.e. $g\in {\cal G}_x^y$. Then the representation $ ({\cal
U},\mathbf{H}_1)$ is called \emph{a subrepresentation of $ ({\cal
U},\mathbf{H})$}. A subrepresentation $ ({\cal U},\mathbf{H}_1)$
is called \emph{proper} if \ $\mathbf{H}_1$ is proper subbundle of
\ $\mathbf{H}$, i.e. $\mathbf{H}_1 \neq \mathbf{H}$ and
$\mathbf{H}_1$ is not null space bundle.
\end{Definition}
\begin{Definition}
A unitary representation $ ({\cal U},\mathbf{H})$ is called
\emph{irreducible} if it has no proper subrepresentations.
\end{Definition}
\begin{Example}
Let $\mathbf{H} = X\times H$ be a trivial Hilbert bundle over $X$
with fiber $H$. For $g\in {\cal G}_x^y, \ x,y\in X$ define an
operator of the representation \\ $U(g): \{x\}\times H\rightarrow
\{y\}\times H$ by \ $U(g)(x,h)= (y,h)$. Such representation
$({\cal U},\mathbf{H})$ of the groupoid $\cal G$ is called
\emph{trivial representation}. A trivial representation is
irreducible if and only if it is one-dimensional, i.e., $dim H =
1$.
\end{Example}
\begin{Example}
Let $H_x = L^2({\cal G}_x, d\lambda_x)$, for $x\in X$, be a
Hilbert space of square $\lambda_x$-integrable functions on ${\cal
G}_x$, and for $g\in {\cal G}_x^y, x,y\in X$ and $f \in H_x$
define $U(g): H_x \rightarrow H_y$ by
\[ (U(g)f)(g_1)= f(g_1\circ g),\]
for $g_1 \in {\cal G}_y$.
\par
A representation $({\cal U},\mathbf{H})$ is called \emph{regular
representation of the groupoid $\cal G$ \cite{Paterson}.}
\end{Example}
Now let us consider the regular representation of a pair
groupoid ${\cal G}_0 = X\times X$. Let $\mu$ be a regular Borel
measure on $X$. Let us define a right Haar system of measures
$\{\mu_x\}_{x\in X}$ on the pair groupoid ${\cal G}_0$, $\mu_x$
being given on ${\cal G}_{0,x} = X \times \{x\}$ by the formula
$\mu_x(f) = \int_{{\cal G}_{0,x}} f(y, x)d\mu(y)$.
\par
Then we have a simple invariance condition:
\[ \int_{{\cal G}_{0,x}} f[(y,x)\circ(x,z)]d\mu(y)=\int_{{\cal
G}_{0,z}} f(y,z)d\mu(y). \]

For each $x \in X$ the Hilbert space $L^2({\cal
G}_{0,x},d\mu_x)$ is obviously isomorphic to $L^2(X)$.
\begin{Example}
The regular representation of a pair groupoid ${\cal G}_0$ in the
Hilbert bundle $\{ L^2({\cal G}_{0,x})\}_{x \in X}$ over $X$ is
given by the following family of operators $\{U_0(g)\}, \ g\in
{\cal G}_{0,x}^y, \  x,y\in X$
\[ [U_0(g) f] (z, y) = f(z, x) \]
where $z\in X$.
\end{Example}
Let us observe that the regular representation of pair groupoid is
equivalent to trivial representation in the trivial Hilbert bundle
$X\times L^2(X)$. \vspace{1cm}
\par
Now, we shall introduce the \emph{quotient groupoid ${\cal
G}/\Gamma$} (cf. \cite{Mackenzie}) and consider its
representations.

Let $\Gamma$ be the isotropy groupoid of a groupoid ${\cal G}$.
Let us define an equivalence relation $\sim$ on ${\cal G}$, for $g, h
\in {\cal G}$,
\[ g \sim h \Leftrightarrow \mathrm{there \ exist} \ \gamma_1 \in
\Gamma \ \mathrm{such \ that} \ (\gamma_1, g) \in {\cal G}^2 \ \mathrm{and} \
\gamma_1 \circ g = h. \]

Let us notice that if $\gamma_2 \in \Gamma_{d(g)}$ then also $g \sim
g \circ \gamma_2$. Indeed, $g \circ \gamma_2 = g \circ \gamma_2
\circ g^{-1} \circ g = \gamma_1 \circ g$ where $\gamma_1 = g \circ
\gamma_2 \circ g^{-1} \in \Gamma_{r(g)}$.
\par
Denote the equivalence class of $g \in {\cal G}$ by $[g]$, and the
set of such equivalence classes by ${\cal G}/\Gamma$. Then we can
introduce the groupoid structure on ${\cal G}/\Gamma$. The
structure maps $\tilde{d}$ and $\tilde{r}$, the multiplication, the
inverse and the identity section $\tilde{\epsilon}$ are given by
$\tilde{d}[g]=d(g)$, $\tilde{r}([g])=r(g)$, $[g] \circ [h]=[g \circ
h]$, for $(g, h) \in {\cal G}^2$, $[g]^{-1} = [g^{-1}]$,
$\tilde{\epsilon}(x) = [\epsilon(x)]$, respectively.
\par
It easy to see that the canonical projection $p: {\cal G}
\rightarrow {\cal G}/\Gamma$ is a homomorphism of (topological)
groupoids (in ${\cal G}/\Gamma$ we choose the quotient topology).
Let us denote by ${\cal G}_0$ the pair groupoid ${\cal G}_0 = X
\times X$ over the base $X$. Recall that in ${\cal G}_0$ we have
$d_0 (x,y) = y, r_0(x,y) = x, (x, y) \circ (y, z) = (x, z), (x,
y)^{-1} = (y, x)$ and $\epsilon_0 (x) = (x, x)$.
\par
We observe that the quotient groupoid ${\cal G}/\Gamma$
coincides with ${\cal G}_0$.

\begin{Proposition}
The map $\Phi: {\cal G}/\Gamma \rightarrow {\cal G}_0$, given by
$\Phi([g]) = (r(g), d(g))$, is an isomorphism of groupoids over
$X$.
\end{Proposition}
\par
\textbf{Proof:} It is clear that, for $g, h \in {\cal G}^2$, one has
$\Phi([g] \circ [h]) = (r(g), d(h)) = (r(g), d(g)) \circ (r(h),
d(h)) = \Phi ([g]) \circ \Phi ([h])$. Also $\Phi ([g]^{-1}) = (d(g),
r(g)) = (\Phi([g]))^{-1}$. Thus $\Phi$ is a groupoid homomorphism.
It is clear that $\Phi$ maps ${\cal G}/\Gamma$ onto ${\cal G}_0$.
Moreover, if $\Phi([g])$ is a unit element in ${\cal G}_0$, i.e.,
$\Phi([g]) = (x, x)$ for an element $x \in X$, then
$d(g)=r(g)=x$, i.e., $g \in \Gamma$ and $[g]$ in a unit in ${\cal
G}/\Gamma$. This means that $\Phi$ is an isomorphism. $\diamond$
\par
Now, let us assume that a representation $({\cal U},\mathbf{H})$
of the groupoid ${\cal G}$ is \emph{$\Gamma-$ invariant}, i.e.,
$U(\gamma\circ g)= U(g)$ for every $g\in {\cal G}, \ \gamma\in
\Gamma, \ (\gamma,g)\in {\cal G}^{2}$. Then it is easily seen that
one has a unitary representation $({\cal
U}_0,\mathbf{H})=(\{U_0([g])\}_ {[g]\in {\cal
G}/\Gamma},\mathbf{H} )$ of the groupoid ${\cal G}/\Gamma $ formed
by the family of operators
\[ U_0([g])= U(g):H_x\rightarrow H_y \]
for every $g$ such that $d(g)=x, \  r(g)=y$.
\begin{Example}
Let $H_x = L_0^2({\cal G}_x)$ be a Hilbert space of
$\Gamma-$invariant and $\mu -$square integrable functions on
${\cal G}_x$ , i.e. such that, for Borel-measurable functions $f$
on ${\cal G}_x$, $f(\gamma\circ g)=f(g)$ for every $g\in {\cal
G}_x, \ \gamma\in \Gamma, \ (\gamma,g)\in {\cal G}^{2}$, and
$\int_{X} |f(g)|^2 d\mu(r(g))<\infty$. It is clear that the space
$L_0^2({\cal G}_x)$ is isomorphic to $L^2(X)$. Define for every
$g\in {\cal G}_x^y, \ x,y\in X$ an operator $U(g): H_x \rightarrow
H_y$ by
\[ (U(g)f)(g_1)= f(g_1\circ g)\]
for $f \in H_x$,  $g_1 \in {\cal G}_y$. In such a manner we obtain
a $\Gamma$-invariant unitary representation of the groupoid $\cal
G $ which is called a \emph{quasi-regular representation}. Let us
observe that the corresponding representation ${\cal U}_0$ of the
quotient groupoid ${\cal G}/\Gamma$ coincides with the regular
representation of the pair groupoid ${\cal G}_0$.
\end{Example}

\section{Induced representations of the groupoid ${\cal G}$.}
In this section, we define the representation of ${\cal G}$
induced by a representation of the isotropy subgroupoid $\Gamma$.
From now on we assume that on the groupoid ${\cal G}$ there exists
a right Haar system $\{\lambda_x\}_{x\in X}$ consistent with
Borel regular measure $\mu$ on $X$.
\par
First, we have to construct an appropriate Hilbert bundle.
\par
Assume that there is given a unitary representation $(\tau,
\mathbf{W})$ of the subgroupoid $\Gamma$. Here $\mathbf{W}$ is a
Hilbert bundle over $X$. Let $W_x$ denote a fiber over $x \in X$
which is a Hilbert space with the scalar product $\langle
\cdot,\cdot \rangle_x$, and let $W = \cup_{x \in X} W_x$ denote
the total space of the bundle $\mathbf{W}$.
\par
Let us define, for every $x \in X$, the space ${\cal W}_x$ of
$W$-valued functions $F$ defined on the set ${\cal G}_x$
satisfying the following four conditions:
\begin{enumerate}
\item
$F(g) \in W_{r(g)}$ for every $g \in {\cal G}_x$,
\item
for every $\mu$-Borel measurable $r$-section $s_x: X \to {\cal
G}_x$ (see Proposition 1) the composition $F\circ s_x$ is a
$\mu$-measurable section of the bundle $\mathbf{W}$,
\item
$F(\gamma \circ g) = \tau(\gamma) F(g)$ for $g \in {\cal G}_x, \
\gamma \in \Gamma_{r(g)}$,
\item
$ \int \langle F(s_x(y)), F(s_x(y)) \rangle_y d\mu(y)<\infty $.
\end{enumerate}
\vspace{1cm}
If we identify two functions $F,F' \in {\cal W}_x$
satisfying
\[\int \langle (F-F')(s_x(y)), (F-F')(s_x(y)) \rangle_y
d\mu(y)=0,\]
we can introduce the scalar product $(\cdot,
\cdot)_x$ in the space ${\cal W}_x$
\[ (F_1, F_2)_x = \int \langle F_1(s_x(y)), F_2(s_x(y)) \rangle_y
d\mu(y)
\]
where $s_x$ is the section determined by Proposition 1, part 3.
\par
The spaces ${\cal W}_x$, $x \in X$, with these scalar products are
Hilbert spaces. It is easily seen that they are isomorphic to the
Hilbert space $L^2(X, \mathbf{W})$ of square-integrables sections
of the bundle $\mathbf{W}$. Now, let us denote $\mathbf{\cal W} =
\{ {\cal W}_x \}_ {x \in X} $. It is a Hilbert bundle over $X$. We
define a unitary representation of the groupoid ${\cal G}$ in the
Hilbert bundle $\mathbf{ {\cal W}}$ in the following way

\begin{Definition}
The representation of the groupoid ${\cal G}$ \emph{induced by the
representation $(\tau, \mathbf{W})$} of the subgroupoid $\Gamma$
is the pair $(U^{\tau},\mathbf{\cal W})$ where, for $g\in {\cal
G}_x^y$, we define $U^{\tau}(g):{\cal W}_x\rightarrow {\cal W}_y$
by
\[ (U^{\tau} (g_0) F) (g) = F (g \circ g_0). \]
\end{Definition}
It is clear that $(U^{\tau}, \mathbf{ {\cal W}})$ is a unitary
groupoid representation.
\par
Sometimes we shall use the notation $U^{\tau} = Ind_{\Gamma}^{\cal
G} (\tau)$.

\section{Systems of imprimitivity.}
For a given Hilbert space $W_0$ we can consider the Hilbert space
$L^2(X, W_0)$ of square integrable $W_0$ - valued functions on
$X$. In the space $L^2(X, W_0)$ one has a representation of the
commutative algebra $L^\infty(X)$ given by the multiplication operators
by the function: $L^\infty(X) \ni f \rightarrow
\pi_0(f) \in B(L^2(X,W_0))$ where, for $z \in X$,
\[ [\pi_0(f) \psi](z) = f(z)\psi(z). \]
We shall call $\pi_0$ \emph{the natural representation of
$L^\infty (X)$} in $L^2(X, W_0)$.
\par
Now, let us consider a representation $U$ of the groupoid ${\cal G}$ in a
Hilbert bundle $\mathbf{H}$ over $X$. We assume that for $\mu$ -
a.e. $x \in X$ there exists a Hilbert space $W_x$ with a scalar
product $\langle \cdot, \cdot \rangle_x$ such that the spaces
$W_x$ are isomorphic with each other. Let us assume that, for $\mu$ -
a.e. $x$, the fiber $H_x$ of the bundle $\mathbf{H}$ is isomorphic
to $L^2(X, W_x)$. We shall simply write  $H_x = L^2(X,W_x)$, and
$U(g): L^2(X, W_x) \rightarrow L^2(X,W_y)$ for $g \in {\cal
G}_x^y$. It is clear that the collection of the spaces $\{ L^2(X,
W_x) \}_{x \in X}$ forms a Hilbert bundle over $X$ which is
isomorphic to the bundle $\mathbf{H}$.

\begin{Definition}
We say that there exists \emph{a system of imprimitivity $({\cal
U},\pi)$} for the representation $({\cal U}, \mathbf{H}) $ of the
groupoid ${\cal G}$ if
\begin{enumerate}
\item
the representation $({\cal U}, \mathbf{H}) $ satisfies the above
assumption ($H_x = L^2(X, W_x)$ for $\mu$ -a.e. $x\in X$),
\item
$\pi=(\pi_x)_{x \in X}$ is the family of natural representations
of the algebra $L^\infty(X)$ in the Hilbert spaces $L^2(X, W_x)$,
\item
for every $f \in L^\infty(X)$, and for $\mu$ - a.e. $x, y \in X$, and
$\nu$ - a.e. $g \in {\cal G}_x^y$
\[ U(g)\pi_x(f)U(g^{-1})=\pi_y(f) .\]
\end{enumerate}
\end{Definition}
\begin{Example}
Let $({\cal U},\mathbf{H})$ be the quasi-regular representation of
the groupoid ${\cal G}$, defined in Example 6. Then there exists a
system of imprimitivity $({\cal U},\pi)$ for ${\cal U}$. Indeed, for
$f\in L^\infty(X), \ \psi \in H_y, \ g\in {\cal G}_x^y, \ h\in {\cal
G}_y$, we have
\[(U(g)\pi_x(f)U(g^{-1}) \psi) (h)= (\pi_x(f)U(g^{-1}) \psi) (h\circ g)= \]
\[ =f(r(h\circ g))(U(g^{-1}) \psi) (h\circ g)= f(r(h))\psi (h) = \pi_y(f)\psi(h).\]
The quasi-regular representation can be understood as induced by a
trivial one-dimensional representation of the subgroupoid
$\Gamma$.

\end{Example}
\par
We are now in a position to state our main theorem (the
imprimitivity theorem for groupoids):

\begin{Theorem}
If, for a representation $({\cal U}, \mathbf{H})$, there exists a
system of imprimitivity $({\cal U}, \pi)$ then the representation
${\cal U}$ is equivalent to the representation ${\cal U}^\tau$
induced by some representation $(\tau, \mathbf{W})$ of the
subgroupoid $\Gamma$.
\end{Theorem}
\par
Let us observe that, for $\gamma \in \Gamma_x = {\cal G}_x^x$,
condition 3 of the Definition 9 of the imprimitivity system
reduces to the following one
\[ U(\gamma) \pi_x(f) U(\gamma^{-1}) = \pi_x (f). \]
Let us denote by ${\cal M}_{0,x}$ the following subalgebra in
$B(H_x) = B(L^2(X,W_x))$: ${\cal M}_{0,x} = \{ U(\gamma): \gamma
\in \Gamma_x \}$.
\par
Then we have, for $\mu$ - a.e. $x \in X$, $\pi_x(L^\infty(X))
\subset {\cal M}_{0,x}'$, where ${\cal M}_{0,x}'$ denotes the
commutant of the algebra ${\cal M}_{0,x}$ in $B(H_x)$.

\begin{Definition}
The system of imprimitivity $({\cal U}, \pi)$ is
\emph{irreducible} if, for $\mu$ - a.e. $x \in X$,
$\pi_x(L^\infty(X)) = {\cal M}_{0,x}'$.
\end{Definition}

\begin{Theorem}
If for a representation $({\cal U}, \mathbf{H})$ there exists an
irreducible system of imprimitivity $({\cal U}, \pi)$ then the
representation ${\cal U}$ is equivalent to the representation ${\cal
U}^\tau$ induced by some irreducible representation $(\tau,
\mathbf{W})$ of the subgroupoid $\Gamma$.
\end{Theorem}

First, let us notice that from the fact that all operators
$U(\gamma), \ \gamma \in \Gamma_x$, commute with all $\pi_x(f), \ f
\in L^\infty(X)$, \ it follows that $U(\gamma)$ are decomposable
(see \cite{Dixmier}, part II, 2.5). This means that, for $\mu - a.e.
\ y \in X$, there exists an operator $U(\gamma)_y \in B(W_x)$ such
that, for $\psi \in L^2(X, W_x)$, $\ \ (U(\gamma)\psi)(y) =
U(\gamma)_y(\psi(y))$. It is easily seen that all $U(\gamma)_y$ are
unitary (cf. \cite{Dixmier}, II.2, Ex.2). We can prove even more.

\begin{Lemma}
If for a representation $(U, \mathbf{H})$ there exists a system of
imprimitivity, then
\begin{enumerate}
\item
there exists a unitary representation $(\tau_x, W_x)$ of the group
$\Gamma_x$ such that $\tau_x(\gamma) = U(\gamma)$ for every
$\gamma \in \Gamma_x$ and $\mu$ - a.e. $x \in X$. (In particular
it means that the function $X \ni y \rightarrow U(\gamma)_y \in
B(H_x)$ is a constant field of operators),
\item
we can define a representation $(\tau, \mathbf{W})$ of the subgroupoid
$\Gamma$ such that, for $\gamma \in \Gamma_x$, $ \ \tau(\gamma) =
\tau_x(\gamma)$,
\item
if the system of imprimitivity $(U, \pi)$ is irreducible then
$(\tau, W)$ is an irreducible representation of \ $\Gamma$, i.e.,
for $\mu$ - a.e. $x \in X$, the representations $(\tau_x, W_x)$ of
the groups $\Gamma_x$ are irreducible.
\end{enumerate}
\end{Lemma}

\textbf{Proof:} Notice that the Hilbert space $L^2(X, W_x)$ is
isomorphic to the tensor product of Hilbert spaces $L^2(X)
\bigotimes W_x$. A decomposable operator in such a space has the
form $[A (\psi \otimes h)](y) = \psi(y) \otimes A_y h$. We have to
show that it is of the form $id_{L^2} \otimes A_0$, where $A_0 \in
B(W_x)$. Let $\{\psi_i\}_{i=1,2,...}$ be an orthonormal basis of
the space $L^2(X)$. Consider the unitary operators $U_{ij}$ in the
space $L^2(X)$ defined by $U_{ij}\psi_j=\psi_i$. If $A$ is
decomposable then $A$ commutes with all operators of the form
$U_{ij} \otimes id_{H_x}$. Then it follows that $A = id_{L^2}
\otimes A_0$ by Lemma 2 (\cite{Dixmier}, section I.2.3). For the
operators $U(\gamma), \gamma \in \Gamma_x$ let us denote by
$\tau_x(\gamma)$ the operators $W_x \rightarrow W_x$ such that
$U(\gamma) = id_{L^2} \otimes \tau_x (\gamma)$. It is clear that
all $\tau_x(\gamma)$ are unitary in $W_x$, and $\tau_x(\gamma_1
\circ \gamma_2) = \tau_x(\gamma_1) \circ \tau_x (\gamma_2)$ for
$\gamma_1, \gamma_2 \in \Gamma$. Thus $\tau_x$ is a unitary
representation of the group $\Gamma_x$ in the Hilbert space $W_x$.
This ends the proof of part 1.
\par
Now the assertion 2 of the Lemma is obvious.
\par
To obtain part 3 it is sufficient to see that the condition of
irreducibility of the imprimitivity system implies that only
operators of the form $\lambda id_{W_x}, (\lambda \in \mathbf{C})$
commute with all $\tau_x(\gamma), \gamma \in \Gamma_x$. But by
Schur's lemma it follows that the representation $\tau_x$ of
$\Gamma_x$ is irreducible.  $\diamond$
\par
The next lemma gives us more properties of the representation $(\tau,
\mathbf{W})$ of the groupoid $\Gamma$ as well as of the
representation $(U, \mathbf{H})$ that has a system of
imprimitivity.

\begin{Lemma}
$ $
\begin{enumerate}
\item
The representations $\tau_x, x \in X$ are equivalent to each
other, as representations of isomorphic groups $\Gamma_x$.
\item
The operators $U(g): H_x \rightarrow H_y$, where $H_x = L^2(X,
W_x)$, $H_y = L^2(Y, W_y)$ for $g \in {\cal G}_x^y$, are
decomposable, i.e., there exist unitary operators $U^0(g): W_x \to
W_y$ such that for $\psi \in L^2(X, W_x)$ and, for $z \in X$,
\[ (U(g) \psi)(z) = (U^0(g)) (\psi (z)). \]
Moreover, the operator $U^0(g): W_x \to W_y$ does not depend of $z
\in X$.
\end{enumerate}
\end{Lemma}
\par
\textbf{Proof:} First we shall prove part 2. First of all, let us
notice that all spaces $W_x$, for $\mu$ - a.e. $x \in X$, are
isomorphic to each other as Hilbert spaces. Denote by $i_x^y: W_x
\to W_y$ the isomorphism and define the unitary map $R_x^y: L^2(X,
W_x) \to L^2(X,W_y)$ by $(R_x^y \psi)(z) = i_x^y(\psi(z))$, $\psi
\in L^2(X,W_x), z \in X$. Consider the composition of unitary maps
$U(g) \circ (R_x^y)^{-1}: L^2(X,W_y) \to L^2(X,W_y)$ where $g \in
{\cal G}_x^y$. By using the property of the imprimitivity system for
$U(g)$, we obtain
\[
U(g) \circ (R_x^y)^{-1} \circ \pi_y(f) = \pi_y(f) \circ U(g) \circ
(R_x^y)^{-1}
\]
for $f \in L^\infty(X)$.
\par
This means that the operator $U(g) \circ (R_x^y)^{-1}$ is
decomposable in $L^2(X,W_y)$. But $(R_x^y)$ is a decomposable map
by definition, therefore $U(g)$ is decomposable as the composition
of decomposable maps. As in the proof of Lemma 1 we conclude that
$U^0(g)$ does not depend of $z \in X$ and is unitary.
\par
To prove part 1 let us first observe that the isotropy groups
$\Gamma_x$ are isomorphic to each other $x \in X$. Indeed, taking
an element $g \in {\cal G}_x^y$ we define the isomorphism
$i:\Gamma_x \to \Gamma_y$ by the formula $i(\gamma) = g \circ
\gamma \circ g^{-1}$ for $\gamma \in \Gamma_x$. Now, we have
$U(i(\gamma)) = id_{L^2} \otimes \tau_y(i(\gamma))$ as in the
proof of Lemma 1. On the other hand, $U(i(\gamma)) = U(g) \circ
U(\gamma) \circ U(g^{-1}) = (id_{L^2} \otimes U^0(g)) \circ
(id_{L^2} \otimes \tau_x(\gamma)) \circ (id_{L^2} \otimes
U^0(g)^{-1}) = id_{L^2} \otimes (U^0(g) \circ \tau_x(\gamma) \circ
U^0(g)^{-1})$. Therefore, we have $\tau_y(i(\gamma)) = U^0(g)
\circ \tau_x(\gamma) \circ U^0(g)^{-1}$, but this means that the
representations $\tau_y$ and $\tau_x$ are equivalent.
\\
\par
Now, we are in a position to give proofs of Theorems 1 and 2.
\par
\textbf{Proof}. Let us consider the spaces $\{ {\cal W}_x \}_{x
\in X}$, introduced in Section 1, connected to the representation
$\tau$ of Lemma 1 and the corresponding induced representation
$U^\tau$. We shall show that the representation $(U, \mathbf{H})$
is equivalent to $(U^\tau, {\cal W})$. We define a family of
isomorphisms of Hilbert spaces $J_x: H_x \to {\cal W}_x$ for $\mu$
- a.e. $x \in X$. Since $H_x=L^2(X,W_x)$, for $\psi \in H_x$, $g
\in {\cal G}_x$, and $r(g)=y$, we put $F(g)=(J_x \psi)(g) =
(U(g)(\psi))(y)$. The definition is correct since by Lemma 2 we
have $(U(g) \psi) (y) = U^0 (g) (\psi(y))$, and $U^0 (g)$ does not
depend of $y \in X$. Since $U(g) \psi \in L^2(X,W_y)$, therefore
$[U(g)(\psi)](y) \in W_y$. Also it is clear that $F(\gamma \circ
g) = \tau (\gamma) (F(g))$ for $\gamma \in \Gamma_y$. To see the
square-integrability let us write
\[
\int\langle F(s_x(y)), F(s_x(y))\rangle_{y} d\mu (y) =
\]
\[
= \int \langle
U^0(s_x(y))(\psi)(y),U^0(s_x(y))(\psi)(y)\rangle_{y} d\mu (y) =
\int \langle \psi(y), \psi(y)\rangle_{y} d\mu(y) =
\]
\[
= \parallel \psi
\parallel_{H_x} < \infty.
\]
This also shows that $J_x$ are unitary maps and are injective. To
see that $J_x$ map onto ${\cal W}_x$, we can give the formula for
$J_x^{-1}$: \ $(J_x^{-1}F)(y)=(U^0(g))^{-1}(F(g))$ where $F \in
{\cal W}_x$ and $g \in {\cal G}_x^y$. Then the right-hand side does
not change if we take other element $g_1 \in {\cal G}_x^y$. Indeed,
since $g_1 = \gamma \circ g$, for an element $\gamma \in \Gamma_y$,
therefore we have $(U^0(\gamma \circ g))^{-1}(F(\gamma \circ g)) =
((U^0(g))^{-1} (\tau(\gamma))^{-1}(\tau(\gamma)) (F(g)) =
(U^0(g))^{-1}(F(g))$. This shows that $J_x$, $x \in X$, are
isomorphisms of Hilbert spaces. Now we can see that $J_x$ are
intertwining maps for the representations $U$ and $U^\tau$, i.e.,
that the following diagram commutes
\[
\begin{CD}
H_x @>{U(g)}>> H_z \\
@V{J_x}VV @VV{J_z}V \\
{\cal W}_x @>>{U^\tau(g)}> {\cal W}_z
\end{CD}
\]
for $\mu$-a.e. $x,z \in X$ and ${\cal \nu}$ - a.e. $g \in {\cal
G}_x^z$. Let $\psi \in H_x$. Then, for $h \in {\cal G}_z^y$, we have
$[(J_z U(g))(\psi)](h)=[(U(h)(U(g))(\psi)](y) = U(h \circ g)(\psi
(y)) = U^0(h \circ g) (\psi (y))$. On the other hand, $ U^{\tau}(g)
J_x (\psi)(h) = [J_x (\psi)](h \circ g) = [U(h \circ g) (\psi)](y)$.
This ends the proof of Theorem 1.
\par
The theorem 2 is now a simple consequence of Theorem 1 and Lemma 1,
part 3.

\section{Representations of the transformation groupoid ${\cal G} = X
\times G, X = K \backslash G$ }

As an introduction to this section we recall the concept of
\emph{induced representation in Mackey sense} (cf \cite{Mackey3},
\cite{Mackey1}, \cite{Taylor}) of a Lie group $G$ by a unitary
representation $(L, V)$ of its closed subgroup $K$ defined in a
Hilbert space $V$.
\par
We assume, for simplicity, that $X = K \backslash G$ has a
$G$-invariant measure $\mu$. We consider ${\cal H}_L$, a Hilbert
space consisting of measurable functions $\phi$ on $G$ with values
in $V$, such that
\[ \phi(h g) = L(h) \phi(g), h \in K, \]
and
\[ \int_X ||\phi([g])||_V^2 d\mu([g]) < \infty \]
where $[g]$ denotes the image of $g$ in $X$ under the projection
$G \to K \backslash G = X$. We introduce the inner product
\[ (\phi_1, \phi_2)_{{\cal H}_L} = \int_X (\phi_1(x), \phi_2(x))_V d\mu(x). \]
Then we define the representation $U^L$ of $G$ on ${\cal H}_L$
given by the formula
\[ (U^L(g) f)(g_0) = f(g_0 g), \ g_0,g \in G, \ f \in {\cal H}_L. \]
It is easily seen that $U^L$ is unitary. The representation $(U^L,
{\cal H}_L)$ is called \emph{induced by the representation $L$ of
$K$}.

\par
Let $G$ be a noncompact Lie group and $K$ its compact subgroup. We
assume that $G$ is unimodular. Then the homogeneous space $X = K
\backslash G$ is a $G-$manifold with right action of the group $G$:
$X \times G \ni (x,g) \mapsto x  g \in X$.
\par
As above, we assume that there exists a $G$-invariant measure $\mu$
on the space $X$, i.e., for $f \in C_c(G)$ we have
\[ \int f(x  g)d\mu(x) = \int f(x)d\mu(x). \]
We shall consider the structure of transformation groupoid on
${\cal G}= X\times G$ (cf. Example 2) and construct a right Haar
system on ${\cal G}$ consistent with the measure $\mu$.
\par
Let us denote
\[ {\cal G}_x = \{ (x,g) \in {\cal G}: g \in G\}, \]
\[ {\cal G}^y = \{ (yg^{-1},g) \in {\cal G}: g \in G\}. \]
Let us also denote the isotropy group ${\cal G}_x^x$ by $\Gamma_x$,
$\Gamma_x = \{(x,k): k \in K_x\}$, where $K_x$ is a subgroup of $G$
of the form $K_x = g_0^{-1} K g_0$ where $g_0 \in G$ is an element
of the coset $x$ ($x = [g_0]$). Indeed, for $k_x \in K_x$ we have $x
k_x = [g_0] g_0^{-1} k g_0 = [ k g_0 ] = x$.

\begin{Lemma}
Let $s_0$ be a Borel section of the principal bundle $G \to K
\backslash G = X$, i. e., $[s_0(x)] = Ks_0(x) = x$. Then
\begin{enumerate}
\item
for every $x \in X$ there exists a section $s_x : X \to {\cal G}_x$
with respect to the map $r$, i.e., $r(s_x(y))=y$,
\item
every element $\mathbf{g} \in {\cal G}_x$ can be represented as
$\mathbf{g} = \mathbf{k}\circ s_x(y)$ where $\mathbf{k} \in
\Gamma_y = {\cal G}_y^y$.
\end{enumerate}
\end{Lemma}
\textbf{Proof:}
\begin{enumerate}
\item
Let $o \in X$ denotes the origin point, i.e., $o = [k]$, $k \in K$.
Then we have $[k s_0(x)] = x$ or, equivalently, $o s_0(x) = x$.
Analogously, $o s_0(y) = y$ for $y \in X$. Thus we can define the
section $s_x: X \to {\cal G}_x$ by the formula $s_x(y) = (x,
s_0(x)^{-1} s_0(y))$. It is clear that $r(s_x(y))=y$.
\item
Let us observe that the product $\mathbf{k}\circ s_x(y)$, where
$\mathbf{k} \in \Gamma_y = {\cal G}_y^y$, is of the form
\[(y,s_0(y)^{-1} k s_0(y)) \circ (x,s_0(x)^{-1} s_0(y)) =
(x,s_0(x)^{-1} k s_0(y)). \]

We have to show that if $\mathbf{g} = (x,g)$ with $g$ such that
$xg=y$, then there exists $k \in K$ such that $g=s_0(x)^{-1} k
s_0(y)$ or, equivalently, $s_0(x) g s_0(y)^{-1} = k$. Now, $o s_0(x) g
s_0(y)^{-1} = x  g  s_0(y)^{-1} = y s_0(y)^{-1}=o$. But the isotropy
group of the origin point $o$ is equal to $K$, what means $s_0(x)
 g s_0(y)^{-1} \in K$. $\diamond$
\end{enumerate}
Now, for a function $f \in C_c({\cal G}_x)$, let us define
$f_x(y,k) = f(x,s_0(x)^{-1} k s_0(y))$, and
\[ \int_{{\cal G}_x} f(\mathbf{g}) d\lambda_x (\mathbf{g}) = \int_X
\int_K f_x(y,k) dk d\mu(y). \]

\begin{Proposition}
The collection $\{\lambda_x\}_{x \in X}$ is a right Haar system on
the groupoid ${\cal G}$ consistent with the measure $\mu$ on $X$.
\end{Proposition}

\textbf{Proof:} We have to show the right-invariance of the system
$\{\lambda_x\}$. Let $z=r(\mathbf{g}_0)$, $\mathbf{g}_0 = (x,g_0)$
such that $x g_0 = z$ and let us compute
\[ \int_{{\cal G}_z} f(\mathbf{g} \circ \mathbf{g}_0) d\lambda_z(\mathbf{g}) =
\int_X \int_K f[(z, s_0(z)^{-1} k s_0(y)) \circ (x,g_0) ] dk d\mu(y)
= \]

\[ \int \int f(x,g_0 s_0 (z)^{-1} k s_0(y)) dk d\mu(y). \]

Let us observe that $os_0(x)g_0 s_0 (z)^{-1} = z s_0(z)^{-1} = o$,
thus $ s_0(x)g_0 s_0 (z)^{-1} \in K$  and $ s_0(x)g_0 s_0 (z)^{-1}
k = k_1 \in K$. But this means that $g_0 s_0 (z)^{-1} k = s_0
(x)^{-1}k_1 $. Thus, continuing the computation, we have
\[ \int \int f(x, g_0 s_0(z)^{-1} k s_0(y)) dk d\mu(y) =
\int \int f(x, s_0(x)^{-1} k_1s_0(y)) dk_1 d\mu(y) = \]

\[ \int_{{\cal G}_x} f(\mathbf{g}) d\lambda_x(g).\hspace{0.5cm} \diamond\]

Now, we shall consider representations of the isotropy subgroupoid
$\Gamma$. As we have seen, $\Gamma = \bigcup_{x \in X} \{x\} \times
K_x$ with $K_x = g^{-1} K g$ and $g \in G$ such that its
coset in $X$ is equal to $x$ ($[g] = x$). We can use $g=s_0(x)$.
\par
Let $(\tau, \mathbf{W})$ be a unitary representation of the groupoid
$\Gamma$ in a Hilbert bundle $\mathbf{W} = \{W_x\}_{x \in X}$.

\begin{Definition}
A representation $(\tau, \mathbf{W})$ is called $X$-consistent if
there exist a unitary representation $(\tau_0,W_0)$ of the group
$K$ and a family of Hilbert space isomorphisms
\[ A_x : W_0 \to W_x, \ x \in X \]
such that, for $\gamma \in \Gamma_x$ of the form $\gamma = (x,
s_0(x)^{-1} k s_0(x))$,
\[ \tau(\gamma) = A_x \tau_0(k)  A_x^{-1}. \]
\end{Definition}

\begin{Proposition}
Let $(U,\mathbf{W})$ be a unitary representation of the
groupoid ${\cal G}$. Then the restriction $(\tau, \mathbf{W})$ to
the subgroupoid $\Gamma$ of the representation $(U,\mathbf{W})$,
given by the formula $\tau(\gamma)=U(\gamma)$ for $\gamma \in
\Gamma$, is a $X$-consistent representation of $\Gamma$.
\end{Proposition}

\textbf{Proof.} We can write $\gamma = (x, s_0(x)^{-1} k s_0(x)) =
(o, s_0(x)) \circ (o, k) \circ (o,s_0(x))^{-1}$ . Then $U(\gamma)
= U((o,s_0(x))) U((o,k)) U((o, s_0(x)))^{-1}$. Let us denote $A_x
= U((o, s_0(x)))$ and $\tau_0(k) = U((o,k)), \ W_0=W_o$. Then it
is clear that $\{A_x\}_{x \in X}$ and $\tau_0$ satisfies the
conditions of $X$-consistent representation.  $\diamond$
\par
In the sequel we shall consider the representation of the groupoid
${\cal G} = X \times G$ induced by $X$ - consistent representation
$(\tau, \mathbf{W})$ of the subgroupoid $\Gamma$, and we shall
establish its connection with the induced representation in the Mackey
sense of the group $G$. We use the notation of section 3. Now
condition 3 of the definition of the space ${\cal W}_x$ assumes the
form
\[ F(\gamma \circ (x,g)) = \tau(\gamma) F(x,g) \]
where $ x,y \in X, \ y=xg, \ g \in G, \ \gamma \in \Gamma_y =
\{y\} \times K_{y} $. Thus we have $\gamma = (y, s_0(y)^{-1} k
s_0(y))$ for an element $k \in K$ .  Then, by the definition of
$X$-consistent representation, we can write
\[ F(\gamma \circ (x,g)) = (A_y \tau_0(k) A_y^{-1}) F(x,g). \]
Let introduce a function $\phi: G \to W_0$ defined by the formula
$\phi(k s_0(y)) = A_{y}^{-1} (F(x,s_0(x)^{-1}k s_0(y)))$. Then the
function $\phi$ has the property $\phi(k g) = \tau_0(k) \phi(g)$.

It is sufficient to check the above formula for $g=s_0(y)$. If $\gamma =
(y, s_0(y)^{-1} k s_0(y))$ then we have $\phi(k g) = A_y^{-1}
(F(\gamma \circ(x,s_0(x)^{-1} s_0(y) ))) = A_y^{-1} (A_y \tau_0(k)
A_y^{-1}) F(x, s_0(y)) = \tau_0 (k) \phi(g)$.
\par
We shall use the notation $(L, W_0)$ for the unitary
representation of the group $K$ in the space $W_0$,  $L = \tau_0$.
Thus we have $\phi(k g) = L (k)\phi(g)$ and we can consider the
Hilbert space ${\cal H}_{L}$ introduced above as well as the
representation $(U^{L}, {\cal H}_{L})$ of the group $G$ induced in
the sense of Mackey by $L$ from the subgroup $K$.
\par
The following theorem establishes a connection of the induced
representation $(\cal U^{\tau}, \cal W)$ of the groupoid ${\cal G}$
with the representation $(U^{L}, {\cal H}_{L})$ of the group $G$.
\par
Denote by $R_g$, $g \in G$, the following operator acting in the
space ${\cal W}_x, x\in X$, $y = xg$,
\[ (R_g F)(x, h) = (A_{xh} A_{xhg}^{-1})
(F(x, h g)). \]
\par
Then we have the family of unitary $G$-representations $(R, {\cal
W}_x), x \in X$. (The unitarity follows from the fact that the measure $\mu$ is
$G$-invariant and the operators $A_{xh}, A_{xhg}$ are Hilbert space
isomorphisms.)

\begin{Theorem}
$ $
\begin {enumerate}
\item
For every $x \in X$ the $G$-representation $(R, {\cal W}_x)$ is
unitarily equivalent to the induced representation $(U^{L}, {\cal
H}_{L})$.
\item
All  representations $(R, {\cal W}_x)$, $x \in X$, are unitarily
equivalent to each other. The equivalence is given by the operators
 $I_x^y: {\cal W}_x \to {\cal W}_y, $
 \[(I_x^y F)(y,s_0(y)^{-1}k s_0(z))= (A_y A_z^{-1})(F(x,s_0(x)^{-1}k s_0(z))),\]
 $x,y \in X$.
\end {enumerate}
\end{Theorem}

\textbf{Proof.}
\begin{enumerate}
\item
We define the linear map $J_x : {\cal W}_x \to {\cal H}_{L}$ by
$(J_x F) (g) = \phi(k s_0(y)) = A_{y}^{-1} (F(x,s_0(x)^{-1}k
s_0(y)))$ where $g=k s_0(y)$. $J_x$ is a linear isomorphism since
$A_{y}$ is an isomorphism and it is easily seen that $J_x$ preserves
scalar products of ${\cal W}_x$ and ${\cal H}_{L}$ and so it is a
Hilbert space isomorphism. To see that it defines an equivalence of
representations, we have to show that, for $g \in G$, the following
diagram is commutative

\[
\begin{CD}
{\cal W}_x @>{R_g}>> {\cal W}_x \\
@V{J_x}VV @VV{J_x}V \\
{\cal H}_{L} @>>{U^{L}(g)}> {\cal H}_{L}
\end{CD}
\]

Let us compute $(U^{L} (g)  J_x)(F)(h) $. It is sufficient to take
$h = s_0(y)$ and to notice that each $g\in G$ can be written in the form
$g = s_0(y)^{-1}k s_0(z)$, for $z \in X, z=yg$ and an element $k\in
K$.
\[ (U^{L} (s_0(y)^{-1}k s_0(z))  J_x)(F)(s_0(y)) = (J_x F)(k s_0(z)) = \]
\[=L(k)A_z^{-1} (F(x,s_0(x)^{-1} s_0(z))). \]
On the other hand
\[ (J_x R_g)(F)(s_0(y)) = A_y^{-1} ((R_gF)(x,s_0(x)^{-1}s_0(y)))= \]
\[=A_y^{-1}(A_y A_z^{-1})(F(x,s_0(x)^{-1}k s_0(z))) = A_z^{-1} \tau(\gamma)(F(x,s_0(x)^{-1} s_0(z)))= \]
\[=A_z^{-1}A_z\tau_0(k)A_z^{-1}(F(x,s_0(x)^{-1} s_0(z))) =L(k)A_z^{-1} (F(x,s_0(x)^{-1} s_0(z))).\]
\item
Now it is a simple observation that $I_x^y = J_y^{-1} J_x $.
$\diamond$
\end{enumerate}

\section{A physical picture. A concept of particle in the representation theory framework}.

In papers (\cite{Full}, \cite{Conceptual}, \cite{Random},
\cite{PysiakTime}) we have studied a model unifying general
relativity and quantum mechanics based on noncommutative geometry.
The principal structure of the model is provided by a transformation groupoid
$\mathbf{G} = E \times G$ where $G$ is the Lorentz group and $E$ is
the principal $G$-bundle over the spacetime $M$ (the  total
space of he bundle is formed by all Lorentz frames at all points of $M$). We
have defined a right action of $G$ on $E$, and the multiplication of
elements of the groupoid is introduced as follows
\[ (pg, g_1) \circ (p, g) = (p, gg_1), \]
$p \in E, g, g_1 \in G$.
\par
The model is reduced to the usual quantum mechanics when an act of
measurement is performed. Then we choose a frame $p \in E$ which
represents a reference frame in which the measurement is done. In
the sequel we consider the situation when we want to observe a
particle from a different reference frame situated at a fixed point
$x \in M$. In such a case, our groupoid reduces to the groupoid $\mathbf{G} =
E_x \times G$ where $E_x$ is the fiber of the bundle $E$ over $x$.
The groupoid $\mathbf{G}$ is transitive and its isotropy
subgroupoid is trivial $\Gamma = E_x \times \{e\}$, where $e \in G$
is the neutral element of the group $G$.
\par

A representation $(\tau, \mathbf{W})$ of $\Gamma$ is realized in a
trivial Hilbert bundle $\mathbf{W} = \{W_p\}_{p \in E_x}$ with $W_p
= \mathbf{C}$ and $\tau(\gamma) = id_{W_p}$ for $\gamma = (p, e)$.
Therefore, the induced representation $(U^{\tau}, W)$ of $\mathbf{G}$ is
simply a regular representation (cf ...). Indeed, for $p \in E$, we have
$\mathbf{G}_p = \{(p, g), g \in G\}$, ${\cal W}_p =
L^2(\mathbf{G}_p) \cong L^2(G)$ and, for $F \in {\cal W}_p$, $h \in
G$,
\[ (U^{\tau} (p, g) F) (pg, h) = F((pg, h) \circ (p, g)) = F(p, gh). \]

Let us notice that this regular representation can serve to define
random operators on the groupoid and then to define the von
Neumann algebra of the groupoid (cf. \cite{Random},
\cite{PysiakTime}, \cite{Conceptual}).
\\
\par

Now we pass to quantum mechanical momentum representation of a
particle with the mass $m$. Having fixed (by an act of
measurement) $p \in E$, we have reduced our initial space to $\{p\}
\times G \cong G$. But we want to consider the energy-momentum
space $H$ of the particle, $H = \{ (p_0, p_1, p_2, p_3) \in
\mathbf{R}^4: p_0^2 - p_1^2 - p_2^2 -p_3^2 = m \}$. We have an
action of the group $G = SL_2 (\mathbf{C})$ on the hyperboloid
$H$ (see \cite{Taylor}).
\par
To describe the action we identify $H$ with the set $\overline{H}$
of hermitian $2 \times 2$-matrices with determinant equal to m,
\[ (p_0, p_1, p_2, p_3) \mapsto \left(
                                  \begin{array}{cc}
                                    p_0 - p_3 & p_2 - ip_1 \\
                                    p_2 + ip_1 & p_0 + p_3 \\
                                  \end{array}
                                \right)
 \]
and we let to act $g \in G$ on $\overline{H}$ to the right in the
following way, $\overline{H} \ni A \mapsto g^* A g \in \overline{H}$.
(It is clear that $det(g^* A g) = det A = m$). Next, we see that the
isotropy group of the element $(p_0,0,0,0)$, $p_0 = \sqrt{m}$ is equal
to $K = SU(2)$. Thus we deduce that the homogoneus space $K
\backslash G$ is diffeomorphic to $H$. We can take the phase space
of a particle of the mass $m$ as the space ${\cal G} = K\backslash
G\times G = H\times G$ and consider the algebraic structure of
transformation groupoid on it.
\\
\par

Let $(\cal U, \cal W)$ be a unitary representation of the
groupoid $\cal G$ in a Hilbert bundle $\cal W$. Assume that there
exists an imprimitivity system $({\cal U},\pi)$ for $(\cal U, \cal
W)$. We say that a particle of mass $m$ is represented by the pair $({\cal
U}, \pi)$. We say that it is an elementary particle if the
imprimitivity system $({\cal U},\pi)$ is irreducible \cite{Mackey1}, \cite{Mackey2}. Equivalently (on the strength of
Theorem 1), we can say that the particle is an induced representation
$({\cal U}^{\tau}, \cal W)$ where $\tau$ is a unitary
representation of the isotropy subgroupoid ${\Gamma}$. In the same
way, we can say that the particle is elementary if the inducing
representation $\tau$ is irreducible and, in turn, this means that the
representation $(L,W_0)$, $L=\tau_0$, of the group $K=SU(2)$
is irreducible. Then the representation $(L,W_0)$ is called
the spin of the particle.


\begin{thebibliography}{14}
\bibitem{Bos}
R. Bos, {\it Continuous representations of groupoids},
arXiv:math/0612639v3 [math.RT].
\bibitem{Brown}
R. Brown, {\it From groups to groupoids}, Bull. London Math. Soc.
19 (1987), 113-134.
\bibitem{Buneci}
M.R. Buneci, {\it Groupoid C*-algebras}, Surveys in Mathematics
and its Applications, ISSN 1842-6298, 1 (2006), 71-98.
\bibitem{SilvaWeinstein}
A. Cannas da Silva and A. Weinstein, {\it Geometric Models for
Noncommutative Algebras}, American Mathematical Society, Berkeley,
(1999).
\bibitem{Dixmier}
J. Dixmier, {\it Von Neumann Algebras}, North Holland Publ. Comp.,
Amsterdam, (1981).
\bibitem{Malicious}
M. Heller, Z. Odrzyg\'o\'zd\'z , L. Pysiak, W. Sasin, {\it
Structure of Malicious  Singularities}, Int. J. Theor. Phys., 42,
(2003),  427 - 441.
\bibitem{Full}
M. Heller, L. Pysiak, W. Sasin, {\it Noncommutative unification of
general relativity and quantum mechanics}, J. Math. Phys., 46
(2005), 122501-15.
\bibitem{Random}
M. Heller, L. Pysiak, W. Sasin, {\it Noncommutative dynamics of
random operators}, Int. J. Theor. Phys., 44 (2005), 619-628.
\bibitem{Conceptual}
M. Heller, L. Pysiak, W. Sasin, {\it Conceptual unification of
gravity and quanta,} Int. J. Theor. Phys., 46 (2007), 2494- 2512.
\bibitem{Aharonov}
M. Heller, Z. Odrzyg\'o\'zd\'z , L. Pysiak, W. Sasin, {\it
Gravitational Aharonov-Bohm Effect}, Int. J. Theor. Phys., 47,
(2008), 2566-2575.
\bibitem{Landsman}
N.P. Landsman, {\it Mathematical Topics between Classical and
Quantum Mechanics}, Springer, New York, (1998).
\bibitem{Mackenzie}
K.C.H. Mackenzie, {\it Lie groupoids and Lie algebroids in
Differential Geometry}, London Math. Society Lecture Notes Series,
124, Cambridge University Press, Cambridge, (1987).
\bibitem{Mackey1}
G.W.Mackey, {\it The relationship between classical mechanics and
quantum mechanics}, Contemporary Math., 214 (1998), 91-109.
\bibitem{Mackey2}
G.W.Mackey, {\it Unitary group representations in physics,
probability and number theory}, Benjamin, Reading, Mass., (1978).
\bibitem{Mackey3}
G.W.Mackey, {\it Induced representations of locally compact groups
I,II }, Acta Math., 55 (1952), 101-139; 58 (1953), 193-221.
\bibitem{Mackey4}
G.W.Mackey, {\it Imprimitivity for representations of locally
compact groups}, Proc. Nat. Acad. Sci. U.S.A. 35 (1949), 537-545.
\bibitem{Packer}
J.A. Packer, {\it Applications of the work of Stone and von
Neumann to the theory of wavelets}, Contemporary Math., 365,
(2004), pp. 253-279.
\bibitem{Paterson}
A.L.T. Paterson, {\it Groupoids, Inverse Semigroups, and Their
Operators Algebras}, Birkhauser, Boston, (1999).
\bibitem{PysiakTime}
L. Pysiak, {\it Time Flow in a Noncommutative Regime}, Internat.
J. Theoret. Phys., 46 (1), (2007), 16 - 30
\bibitem{Pysiak}
L. Pysiak, {\it Groupoids, their representations and imprimitivity
systems}, Demonstratio Mathematica 37 (2004), 661-670.
\bibitem{Renault}
J.N. Renault, {\it A groupoid approach to C*-algebras}, Lecture
Notes in Math. 793, Springer-Verlag, New York, (1980).
\bibitem{Taylor}
M.E. Taylor, {\it Noncommutative Harmonic Analysis}, A.M.S.,
Providence, (1986).
\bibitem{Weinstein}
A. Weinstein, {\it Groupoids: unifying internal and external
geometry}, Contemporary Math. 282, (2001), 1-19.
\bibitem{Westman}
J. Westman, {\it Harmonic analysis on groupoids}, Pacific J.
Math., 27 (1968), 621-632.

\end{thebibliography}
\end{document}